\documentclass[12pt]{article}
\usepackage{amsmath, amssymb}
\usepackage{graphicx}
\usepackage{natbib}
\usepackage{url} 
\usepackage[ruled,vlined]{algorithm2e}
\usepackage{booktabs, breakcites}

\usepackage{xcolor}

\newcommand{\blind}{0}

\renewcommand{\baselinestretch}{1.5}

\begin{document}

\def\spacingset#1{\renewcommand{\baselinestretch}%
{#1}\small\normalsize} \spacingset{1}

\newtheorem{example}{Example}


\if0\blind
{
  \title{\bf On the relationship between \\ Uhlig extended and beta-Bartlett processes}
  \author{V\'ictor Pe\~na  \hspace{.2cm}\\
    Baruch College, The City University of New York \\
    and \\
    Kaoru Irie \\
    Faculty of Economics, University of Tokyo}
  \maketitle
} \fi

\if1\blind
{
  \bigskip
  \bigskip
  \bigskip
  \begin{center}
    {\LARGE\bf Title}
\end{center}
  \medskip
} \fi

\bigskip
\begin{abstract}
Stochastic volatility processes are used in multivariate time-series analysis to track time-varying patterns in covariance matrices. Uhlig extended and beta-Bartlett processes are especially convenient for analyzing high-dimensional time-series because they are conjugate with Wishart likelihoods. In this article, we show that Uhlig extended and beta-Bartlett are closely related, but not equivalent: their hyperparameters can be matched so that they have the same forward-filtered posteriors and one-step ahead forecasts, but different joint (smoothed) posterior distributions. Under this circumstance, Bayes factors can't discriminate the models and alternative approaches to model comparison are needed. We illustrate these issues in a retrospective analysis of volatilities of returns of foreign exchange rates. Additionally, we provide a backward sampling algorithm for the beta-Bartlett process, for which retrospective analysis had not been developed. 
\end{abstract}

\noindent%
\vfill

\spacingset{1.5} 

\section{Introduction}

Time-series with time-varying dependence structures arise naturally in finance, neuroimaging, and online marketing. In these applications, stochastic volatility processes are necessary for successful forecasting and decision making. From a Bayesian perspective, \cite{west2020bayesian} shows that models with conjugate sequential updates are particularly attractive for analyzing high-dimensional time-series, since implementing richly-parametrized models that require Markov chain Monte Carlo methods for posterior inference (e.g., \cite{aguilar2000bayesian}, \cite{nakajima2012dynamic}) may not be computationally feasible. Alternatively, there is a large literature on non-Bayesian methods for modeling high-dimensional time-varying covariance matrices; see the literature review in \cite{windle} or the methods reviewed in \cite{bauwens2006multivariate}.

Two classes of stochastic volatility processes that are conjugate with Wishart likelihoods coexist in the literature: matrix-beta processes, which build upon \cite{uhlig}, and beta-Bartlett processes, which were first used in \cite{betabartlett}. To this date, the most flexible matrix-beta process is the Uhlig extended process \citep{windle}, which is the one we consider herein. Both approaches can be used to model high-dimensional time-series: for example, \cite{casarin2014comment} analyzed a $199$-dimensional time-series with the Uhlig extended process. 
 
Our main contributions are (1) studying the relationship between Uhlig extended and beta-Bartlett processes (Section~\ref{sec:relationship}) and (2) providing the first backward sampler in the literature for beta-Bartlett processes (Section~\ref{sec:BS}). We compare the models in a simple, 3-dimensional foreign exchange rates illustration in Section~\ref{sec:FX}. We end the article with conclusions in Section~\ref{sec:conclusions}. 

\section{Notation and Bartlett decomposition}

 We use the notation $1:n$ to denote the set $\{1, 2, \, ... \, , n\}$. Following \cite{pradowest}, we use the notation $\mathcal{D}_t$ for our ``information set'' at time $t$. Before we observe any data, our prior knowledge is denoted $\mathcal{D}_0$. At time $t \in 1:T,$ the information set is $\mathcal{D}_t = \{\mathcal{D}_0, y_{1:t} \}$. We denote $q$-dimensional normal random variables with mean $\mu$ and covariance matrix $\Sigma$ as $N_q(\mu, \Sigma)$, chi-squared random variables with $k > 0$ degrees of freedom as $\chi^2_k$, and Beta random variables with two shape parameters, $a > 0$ and $b > 0$, as Beta($a,b$). The less common Wishart$_q(k, A)$ and MatrixBeta$_q(n/2, k/2)$ distributions are as defined in \cite{windle}. For extrema, we use the notation $a \wedge b  = \min(a, b)$ and $a \vee b = \max(a,b)$. Finally, we use the notation $\text{uchol}(\cdot)$ for the function that returns the upper-triangular Cholesky factor of a symmetric positive-definite matrix. 

The models we study rely heavily on the Bartlett decomposition of Wishart-distributed matrices, which we now review. Let $W \sim \text{Wishart}_q(k, A)$ be a $q \times q$ random matrix with $k > 0 $ and symmetric positive-definite $A$. Its Bartlett decomposition is $W = (U P)' U P$, where $P = \text{uchol}(A)$ and  $U = (u_{ij})_{i,j \in 1:q}$ is a upper-triangular matrix with entries $u_{ij} \stackrel{\text{iid}}{\sim} N_1(0,1)$ for $i < j$, which are independent of $u_{ii}^2  \stackrel{\text{iid}}{\sim} \chi^2_{k - i +1}$ for $i, j \in 1:q$. 

\section{Uhlig extended and beta-Bartlett processes} \label{sec:relationship}

In this section, we define the Uhlig extended and beta-Bartlett processes and explore their relationship. We assume the readers are familiar with them and refer to \cite{windle} and \cite{betabartlett} for further details.

\cite{windle} extend a model that was originally proposed in \cite{uhlig}. Given $(q \times q)$-dimensional symmetric positive-definite matrices $\{y_t\}_{t \in 1:T}$, the model can be written as
\begin{equation} \label{eq:UE}
y_t \mid \Phi^U_t \stackrel{\text{ind.}}{\sim} \text{Wishart}_q(k, (k \Phi^U_t)^{-1}), \, \, 
\Phi^U_t = ( U^U_{t-1} P^U_{t-1} )' \, \Psi_t \, U^U_{t-1} P^U_{t-1}/\lambda,
\end{equation}
where $\Psi_t \sim \text{MatrixBeta}_q(n/2, k/2)$ and $U^U_{t-1}$ and $P^U_{t-1}$ come from the Bartlett decomposition $\Phi_{t-1}^U = (U_{t-1}^U P_{t-1}^U)' U_{t-1}^U P_{t-1}^U$. The model is completed with the prior $\Phi_{0}^U \mid \mathcal{D}_0 \sim \text{Wishart}_q(n+k, (kD_0^U)^{-1}).$ The hyperparameters are $0 < \lambda < 1$, $n > q-1$, and $k$, which is either a positive integer less than $q$ or a real number greater than $q -1$. We refer to the process on $\{\Phi^U_t\}_{t \in 0:T}$ implied by the model above as to the Uhlig extended (UE) process. The prior distributions and forward-filtered posteriors, as derived in \cite{windle}, are given in Table~\ref{tab:summary}. 

In contrast, the beta-Bartlett (BB) stochastic volatility process \citep{betabartlett} can be written as
\begin{equation} \label{eq:BB}
y_t \mid \Phi^B_t \stackrel{\text{ind.}}{\sim} \text{Wishart}_q(k, (k \Phi^B_t)^{-1}), \, \, \, \Phi^B_t = (\tilde{U}_{t} P^B_{t-1})' \tilde{U}_{t} P^B_{t-1}/b,
\end{equation}
where $P^B_{t-1}$ is defined via the Bartlett decomposition $\Phi^B_{t-1} = ( {U}^B_{t-1} P^B_{t-1})' {U}^B_{t-1} P^B_{t-1} $ and $\tilde{U}_t = (\tilde{u}_{ij,t})_{i,j\in 1:q}$ is constructed by modifying the diagonal elements of ${U}^B_{t-1} = (u_{ij,t-1}^B)_{i,j\in 1:q}$ as explained in Table~\ref{tab:summary}. The hyperparameters of the model are $k > 0$, $0 < \beta < 1$, $0 < b < 1$, and $k_0 > 0$, which appears in the prior $\Phi_0 \mid \mathcal{D}_0 \sim \text{Wishart}_q(k_0, (k D^B_0)^{-1})$ for symmetric positive-definite $D^B_0$. We refer to the process defined on $\{\Phi^B_t\}_{t \in 0:T}$ as to the beta-Bartlett (BB) process. The prior distributions and forward-filtered posteriors with this model can be found in Table~\ref{tab:summary}. 

\spacingset{1}

\begin{table}[htbp]
	\caption{Comparison of Uhlig extended and beta-Bartlett.} \label{tab:summary}
	\resizebox{\columnwidth}{!}{%
			\begin{tabular}{ccc}
				\toprule
				 & Uhlig extended & beta-Bartlett \\
				\midrule \addlinespace[0.5em]
				Likelihood & $y_t \mid \Phi^U_t \stackrel{\text{ind.}}{\sim} \text{Wishart}_q(k, (k \Phi_t^U)^{-1})$ & $y_t \mid \Phi^B_t \stackrel{\text{ind.}}{\sim} \text{Wishart}_q(k, (k \Phi^B_t)^{-1})$ \\ \addlinespace[0.5em]
				State evol. & $\Phi^U_t =  (U^U_{t-1} P^U_{t-1})' \, \Psi_t \, U^U_{t-1} P^U_{t-1}/\lambda$ & $\Phi^B_t = ( \tilde{U}_{t} P^B_{t-1})'  \tilde{U}_{t} P^B_{t-1}/b$ \\  \addlinespace[0.5em]
			 	& $\Phi^U_{t-1} = ( U^U_{t-1} P^U_{t-1})' U^U_{t-1} P^U_{t-1}$ & $\Phi^B_{t-1} = (U^B_{t-1} P^B_{t-1})'  U^B_{t-1} P^B_{t-1}$ \\ \addlinespace[0.5em]
				Error & $\Psi_t \sim \text{MatrixBeta}_q(n/2, k/2)$ & $\tilde{u}_{ij, t} = u^B_{ij, t-1}$ \ $(i\not= j)$ \\ \addlinespace[0.5em]
				& & $(\tilde{u}_{ii, t})^2 = \eta_{i, t} (u^B_{ii,t-1})^2$ \\ \addlinespace[0.5em]
				& & $\eta_{i,t} \stackrel{\text{ind}}{\sim} \text{Beta}( ( \beta k_{t-1}-i+1)/2, (1-\beta) k_{t-1}/2)$ \\ \midrule \addlinespace[0.5em]
				Prior at $t$ & $\Phi^U_t \mid \mathcal{D}_{t-1} \sim \text{Wishart}_q(n, (k \lambda {{D^U_{t-1}}})^{-1} )$ & $\Phi^B_t \mid \mathcal{D}_{t-1} \sim \text{Wishart}_q( \beta k_{t-1}, (k b D^B_{t-1})^{-1} )$ \\  \addlinespace[0.5em]
				Post. at $t$ & $\Phi^U_t \mid \mathcal{D}_t \sim \text{Wishart}_q(n+k, (k D^U_{t})^{-1})$ & $\Phi^B_t \mid \mathcal{D}_t \sim \text{Wishart}_q( k_t, (k D^B_{t})^{-1})$ \\  \addlinespace[0.5em]
				& ${D}^U_{t} = \lambda D^U_{t-1} + y_t$ & ${D}^B_{t} = b D^B_{t-1} + y_t$ and
				$k_t = \beta k_{t-1}  + k$ \\  \addlinespace[0.5em]
				\bottomrule
			\end{tabular}
	} 
\end{table}

\spacingset{1.5}

The priors, forward-filtered posteriors, and one-step ahead forecast distributions of the models defined in Equations (\ref{eq:UE}) and (\ref{eq:BB}) coincide under the condition
\begin{equation}
    k_0 = n+k,\qquad \beta = n/(n+k),\qquad b = \lambda\qquad \text{and}\qquad D^B_0 = D^U_0. \label{hpmain}
\end{equation}
The change of variables is bijective, so if the hyperparameters are set by maximizing the marginal likelihoods of the models, the condition is satisfied.  


However, UE and BB aren't equivalent under Equation~(\ref{hpmain}) because $\Phi^U_{t} \mid \Phi^U_{t-1}, \mathcal{D}_{t-1}$ isn't equal in distribution to $\Phi^B_{t} \mid \Phi^B_{t-1}, \mathcal{D}_{t-1}$, {{which we prove by showing that $E(\Phi^U_{t} \mid \Phi^U_{t-1}, \mathcal{D}_{t-1}) \neq E(\Phi^B_{t} \mid \Phi^B_{t-1}, \mathcal{D}_{t-1})$.}} Assume
Equation~(\ref{hpmain}) holds and $\Phi^U_{t-1} = \Phi^B_{t-1} = \Phi_{t-1} = (U_{t-1} P_{t-1})'U_{t-1} P_{t-1}$, with $U_{t-1} = (u_{ij,t-1})_{i,j\in 1:q}$.  {{While $E(\Phi_t^U \mid \Phi_{t-1}, \mathcal{D}_{t-1})$ can be found immediately using the state evolution described in Table~\ref{tab:summary} and Theorem 3.2 in \cite{konno1988exact}, the derivation of $E(\Phi_t^B \mid \Phi_{t-1}, \mathcal{D}_{t-1})$ is more elaborate (see Appendix~\ref{derivations}). The conditional expectation of the difference is}} 
\begin{align*}
E(\Phi_t^U - \Phi_t^B \mid \Phi_{t-1}, \mathcal{D}_{t-1}) &=  \frac{1}{\lambda} P_{t-1}' \left[  \frac{n}{n+k}  U_{t-1}'U_{t-1} -  E(\tilde{U}_t'\tilde{U}_t \mid \Phi_{t-1}, \mathcal{D}_{t-1}) \right] P_{t-1} \\
E[(\tilde{U}_t'\tilde{U}_t)_{ij} \mid  \Phi_{t-1}, \mathcal{D}_{t-1}] &=  \sum_{l=1}^{i \wedge j -1} u_{li , t-1} u_{lj, t-1} + \delta_{ij} \, \frac{(n-i+1) u_{ii, t-1}^2}{n-i+1+k}  + (1-\delta_{ij}) \,  g(i, j)  \\
g(i, j) &= \frac{\Gamma\left(\frac{n-i \wedge j +2}{2}\right)\Gamma\left(\frac{n-i \wedge j +k+1}{2}\right)}{\Gamma\left(\frac{n-i \wedge j +1}{2}\right) \Gamma\left(\frac{n-i \wedge j +k+2}{2} \right)} \, u_{i \wedge j, i \wedge j, t-1} u_{i \wedge j, i \vee j, t-1}, 
\end{align*}
where $(\tilde{U}_t'\tilde{U}_t)_{ij}$ is the $(i,j)^{th}$ element of $\tilde{U}_t'\tilde{U}_t$, $\delta _{ij}=1$ if $i=j$ and $\delta _{ij}=0$ otherwise. 
In general, $E[(\tilde{U}_t'\tilde{U}_t)_{ij} \mid \Phi_{t-1}, \mathcal{D}_{t-1}]$ and $n U_{t-1}' U_{t-1}/{(n+k)}$ aren't equal: if $U_{t-1}$ is diagonal, $E[(\tilde{U}_t'\tilde{U}_t)_{ii} \mid \Phi_{t-1}, \mathcal{D}_{t-1}] = (n-i+1) u_{ii,t-1}^2/(n-i+1+k) \neq n u_{ii, t-1}^2/(n+k)$ for $i \in 2:q$. 

This distinction affects the smoothed posterior distributions $ \Phi^U_{0:T} \mid \mathcal{D}_T$ and $ \Phi^B_{0:T} \mid \mathcal{D}_T$. Dropping process superscripts, the posterior distribution can be factorized as 
\begin{equation} \label{eq:post}
p(\Phi_{0:T} \mid \mathcal{D}_T ) = p( \Phi_T \mid \mathcal{D}_T)  \,  \prod_{t = 1}^{T-1} p(\Phi_{t} \mid \Phi_{t+1}, \mathcal{D}_{t}).
\end{equation} 
The conditionals $p(\Phi_{t} \mid \Phi_{t+1}, \mathcal{D}_{t})$ of UE and BB are different. For UE, we have $\Phi _t = \lambda \Phi _{t+1} + Z_t$, where $Z_t\sim \text{Wishart}_q( k, (k D^B_{t})^{-1})$; for BB, see Section~\ref{sec:BS}. 
{{Below, we compare the conditional expectations and variances of the conditionals in an example to build intuition.}}


\begin{example} \normalfont
Assume Equation~(\ref{hpmain}) holds, let $k = 1$,  $P_t = \text{uchol}( (k D_t)^{-1} )$ be the identity matrix, and $\Upsilon = \text{uchol}(\Phi_{t+1}) = (\upsilon_{ij})_{i,j \in 1:q}$. Then,   
\begin{align*}
E[ (\Phi^U_t)_{ij} \mid \Phi_{t+1}, \mathcal{D}_t ] &=  \lambda \sum_{l=1}^{i \wedge j} \upsilon_{li} \upsilon_{lj} + \delta_{ij}; \, \, \, V[ (\Phi^U_{t})_{ij} \mid \Phi_{t+1}, \mathcal{D}_t] = 1 + \delta_{ij},
\end{align*}
where $(\Phi^U_t)_{ij}$ is the $(i,j)^{th}$ element of $(\Phi^U_t)$. Similarly, for BB:
\begin{align*}
E[ (\Phi^B_t)_{ij} \mid \Phi_{t+1}, \mathcal{D}_t] &= \lambda \sum_{l=1}^{i \wedge j - 1} \upsilon_{li} \upsilon_{lj} + \delta_{ij}(\lambda \upsilon_{ii}^2 + 1) + (1 - \delta_{ij}) \, h(i,j)  \\
V[ (\Phi_t)^B_{ij} \mid \Phi_{t+1}, \mathcal{D}_t] &= 2 \delta_{ij} + (1 - \delta_{ij}) [\lambda^2 \, \upsilon_{i \wedge j, i \vee j}^2 \, \upsilon^2_{i \wedge j, i \wedge j} + \lambda \upsilon_{i \wedge j, i \vee j}^2 - h(i,j)^2 ],
\end{align*}
with $h(i,j) = \sqrt{2 \lambda} 
\, \upsilon_{i \wedge j, i \vee j} \, U(-1/2, 0, \lambda \upsilon_{i \wedge j, i \wedge j}^2/2) $, where $U(a, b, z)$ is Tricomi's confluent hypergeometric function (see e.g. \cite{abramowitz1988handbook}). The expressions  for the diagonal elements coincide but that need not be the case for the off-diagonal elements: $\upsilon_{ij} > 0$ implies $E[ (\Phi_t)^B_{ij} \mid \Phi_{t+1}, \mathcal{D}_t] >  E[ (\Phi_t)^U_{ij} \mid \Phi_{t+1}, \mathcal{D}_t]$ and $\lambda \upsilon_{ij}^2 < 1$ implies  $V[ (\Phi_t)^B_{ij} \mid \Phi_{t+1}, \mathcal{D}_t] < V[(\Phi_t)^U_{ij} \mid \Phi_{t+1}, \mathcal{D}_t]$. {{Derivations of the formulas in this example can be found in Appendix~\ref{derivations}}.}
\end{example}

{{If Equation~(\ref{hpmain}) is satisfied, the marginal likelihoods of UE and BB are equal and Bayes factors cannot be used to compare them.  However, $\Phi_{0:T}^U \mid \mathcal{D}_T$ and $\Phi_{0:T}^B \mid \mathcal{D}_T$ can be substantially different in practice, as we see in Section~\ref{sec:FX}. 

Instead of Bayes factors, we can use posterior likelihood ratios  \citep{aitkin1991posterior} and posterior predictive checks  \citep{gelman1996posterior} to compare the models. Both of these approaches can be implemented given posterior draws, but they have been criticized for, among other reasons, using the data twice \citep{gelman2013inherent}. Alternatively, \cite{kamary2014testing} propose comparing models via mixtures, which here amounts to fitting 
$$
y_t \mid \alpha, \Phi^U_t, \Phi^B_t \stackrel{\text{ind.}}{\sim} \alpha \, \text{Wishart}_q(k, (k \Phi^{U}_t)^{-1}) + (1-\alpha) \, \text{Wishart}_q(k, (k\Phi^{B}_t)^{-1}),  
$$
where $\alpha \sim \text{Beta}(a_0,b_0)$, and studying the posterior distribution of the mixture weight $\alpha$. We implement all of these approaches in an application in Section~\ref{sec:FX}.}}

\section{Backward sampling for beta-Bartlett processes} \label{sec:BS}

Forward-filtered posteriors and forecast distributions for BB were derived in \cite{quintana2010futures}, but a backward sampler was not developed. Here, we present a novel sampler which uses the factorization of $\Phi_{0:T} \mid \mathcal{D}_T$ in Section~\ref{sec:relationship} and consists in drawing $\Phi_T^\ast \sim \text{Wishart}_q(k_T, (k D_T^{B})^{-1})$ and iteratively sampling $\Phi_t^\ast \sim \Phi_t \mid \Phi^\ast_{t+1}, \mathcal{D}_t$. 

Given $\Phi_{t+1}$ and $\mathcal{D}_t$, consider the decomposition $\Phi _{t+1} = (\tilde{U}_{t+1}^{\ast} P_t)' \tilde{U}_{t+1}^{\ast} P_t / b.$ That is,  $\tilde{U}_{t+1}^{\ast} = \text{uchol} ( b (P_t^{-1})'  \Phi _{t+1}P_t^{-1} )$. Then, we can generate $U_t^{\ast} = (u^\ast_{ij,t})_{i,j \in 1:q}$ as follows. The off-diagonal elements are $u^{\ast}_{ij,t} = \tilde{u}^{\ast}_{ij,t+1}$ for $i<j$ and the diagonal elements are $(u^{\ast}_{ii,t})^2 = (\tilde{u}^{\ast}_{ii,t+1})^2 + \theta _{it}$, where  $\theta _{it} \ \stackrel{\text{iid}}{\sim} \ \chi^2_{ (1- \beta) k_t}$. Finally, we can set $\Phi _t = (U_t^{\ast} P_t)'U_t^{\ast} P_t.$

The expression for the conditional of $(u^\ast_{ii,t})^2$ given $(\tilde{u}^\ast_{ii, t+1})^2$ can be justified using standard results for the univariate gamma-beta discount model (see e.g. Exercise 4 in Section 4.6. of \cite{pradowest}). To relate $U_t^{\ast}$ to $\tilde{U}_t^{\ast}$, observe that $\Phi _t = (U_t^{\ast} P_t)'U_t^{\ast} P_t =  (\tilde{U}_t^{\ast} P_{t-1} )' \tilde{U}_t^{\ast} P_{t-1}/b.$ Therefore, 
\begin{equation*}
	\tilde{U}_t^{\ast}  = \text{uchol} ( b (P_{t-1}')^{-1}\Phi _tP_{t-1}^{-1} ) = \text{uchol} ( b  (U_t^{\ast}P_t P_{t-1}^{-1})'U_t^{\ast} P_t   P_{t-1}^{-1} ) =  \sqrt{b} \, U_t^{\ast} P_t P_{t-1}^{-1}.
\end{equation*}

The matrix $P_{t-1}$ is upper-triangular, so it can be inverted at quadratic computational cost using back-substitution. The backward sampler for the UE process requires simulating Wishart random matrices for all $t \in 0:T$. On the other hand, the BB process only requires sampling a Wishart random matrix for $t = T$, and, for $t \in 0:(T-1)$, it requires $q$ chi-squared random variates. Explicit pseudocode for the backward sampler can be found in Algorithm~\ref{fig:BS}. The sampler can be adapted for general multivariate dynamic linear models with BB stochastic volatilties (as in Section 10.4.8 in \cite{pradowest}).

\spacingset{1} 

\begin{algorithm}[h!]
\footnotesize
\KwIn{$b$, $\beta$, $k_t$, and $P_{t} = \text{uchol}( (k D^B_t)^{-1} )$ from $\Phi^B_t \mid \mathcal{D}_t \sim \text{Wishart}(k_t, (k D^B_t)^{-1})$, $t \in 0:T$.}
\KwOut{$ \Phi_{0:T}^{\ast B} \sim \Phi_{0:T}^{B} \mid \mathcal{D}_T$. }
 $U^\ast_T = (u^\ast_{T,ij})_{i,j \in 1:q}$; $u^\ast_{T, ii} \stackrel{\text{ind.}}{\sim}  \chi^2_{k_T-i+1}$; $u^\ast_{T, ij} \stackrel{\text{iid}}{\sim} N_1(0,1)$, $i \in 1:q$ and $i < j \leq q$\;
 $\Phi^{\ast B}_T =  (U^\ast_T P_T)' U^\ast_T P_T$\;
 $\tilde{U}^\ast_{T} = \sqrt{b} \, U^\ast_{T} P_{T} P_{T-1}^{-1}$\;

 \For{descending $t \in T:1$}
{
  $\theta_{(t-1),i} \stackrel{\text{ind}}{\sim} \chi^2_{(1-\beta)k_{t-1}} $ for $i \in 1:q$\;
  $U^\ast_{t-1} = (u^\ast_{(t-1),ij})_{i,j \in 1:q}$\;
  $u^{\ast}_{(t-1), ii} = \sqrt{(\tilde{u}^\ast_{t,ii})^2 + \theta_{(t-1),i}}$,  $i \in 1:q$; $u^\ast_{(t-1),ij} = \tilde{u}^{\ast}_{t,ij}$, $i \in 1:q$ and $i < j \leq q$\;
  $\Phi^{\ast B}_{t-1} =   (U^{\ast}_{t-1} P_{t-1})' U^{\ast}_{t-1}  P_{t-1}$\;
  \If {$t \ge 2$}
 	{$\tilde{U}^\ast_{t-1} = \sqrt{b} \, U^\ast_{t-1} P_{t-1} P_{t-2}^{-1}$\;}
 }
 \Return {$\Phi^{\ast B}_{0:T}$}\;
    \caption{\bf Backward sampler for $\Phi^B_{0:T} \mid \mathcal{D}_T$ }
    \label{fig:BS}
\end{algorithm}
\spacingset{1.5}
 
\section{Illustration: foreign exchange rates} \label{sec:FX}

We perform a retrospective analysis of volatilities of daily returns of exchange rates of three currencies measured in US dollars: euros (EUR), British pounds (GBP), and Canadian dollars (CAD), observed from January 2008 to October 2010 ($T=739$).  The vector of returns $r_t$ can be turned into a rank-1 symmetric matrix by computing $y_t = r_t r_t'$.  Our observational model is $y_t \mid \Phi_t \stackrel{\text{ind.}}{\sim} \text{Wishart}_q(1,\Phi_t^{-1})$, which is equivalent to modeling $r_t \mid \Phi_t \stackrel{\text{ind.}}{\sim} N_q( 0_q, \Phi_t^{-1}).$

The estimate of the volatility matrix at the starting point, $D_0$, is computed as the sample average of the data in 2007. The other hyperparameters are obtained by maximizing the marginal likelihood of the model and are $n=5$ and $\lambda =0.799$. 

\begin{figure}[!htbp]
	\centering
	\includegraphics[width=4.5in]{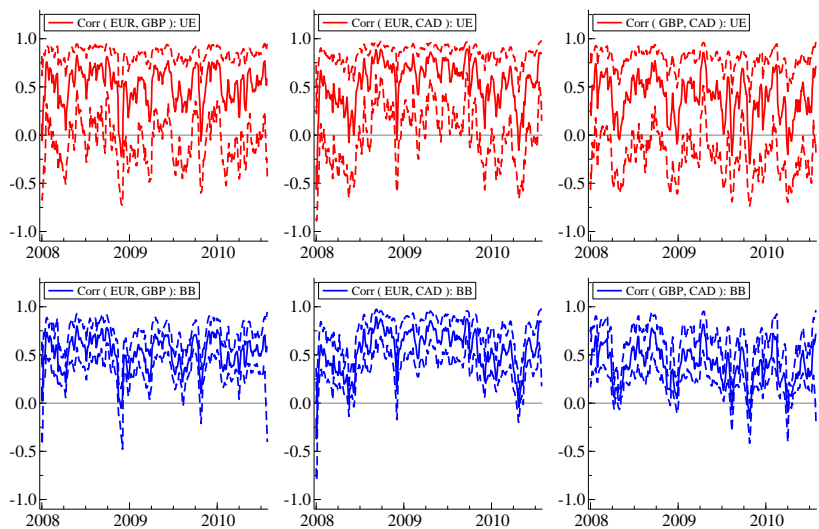}
	\caption{Posterior medians and 95\% credible intervals of correlations computed from sampled $(\Phi _t^U)^{-1}$ and $(\Phi _t^B)^{-1}$ for the UE (top row) and BB (bottom).} \label{fig:corr}
\end{figure}

While UE and BB yield similar point estimates, they are markedly different in their retrospective uncertainty quantification. This difference is apparent in the posterior  correlations displayed in Figure~\ref{fig:corr}. 

The logarithm of the posterior likelihood ratio of the UE model to the BB model is
$$
\overline{\ell}_{UB} = \log\left\{E[ L(\Phi^U_{0:T}) | \mathcal{D}_T]/E[ L(\Phi^B_{0:T}) | \mathcal{D}_T] \right\}, \, \, \, L(\Phi_{0:T}) = \prod_{t=1}^T N_q(r_t \mid 0_q, \Phi_t^{-1}),
$$
where the expectations are computed by posterior samples $\Phi^U_{0:T} \mid \mathcal{D}_T$ and $\Phi^B_{0:T} \mid \mathcal{D}_T$. In this application, $\overline{\ell}_{UB} = -35.230$, which clearly favors the BB model. 

We implement a mixture model to compare the models, as proposed in \cite{kamary2014testing}, by running a missing-data augmented Gibbs sampler (see Appendix) for $N = 10^4$ iterations. A mixture weight $\alpha$ close to $0$ favors the BB model, whereas a mixture weight near $1$ favors UE. Starting with $\alpha \sim \text{Beta}(1,1)$, we estimate $E(\alpha \mid \mathcal{D}_T) = 0.498$ with an estimated standard error of $0.00026$ and $P( \alpha < 0.5 \mid \mathcal{D}_T) = 0.533$ with an estimated standard error of $0.0055165$ (we used batch means estimators for the standard errors; see e.g. \cite{geyer1992practical}). The simulation error is small enough to be confident that the mixture model prefers BB, but the evidence isn't nearly as overwhelming as it is with the logarithm of the posterior likelihood ratio. 


We also compare the models with posterior predictive checks. We found the length of 95\% posterior predictive intervals at each time point, while monitoring their coverage. The BB model has smaller interval lengths at essentially no cost; the cumulative coverage rates of the intervals are higher than 95\% and comparable to those of the UE model. 

{{In Appendix~\ref{appFX}, we compare UE and BB to a Bayesian dynamic factor model based on \cite{aguilar2000bayesian}. The posterior correlations estimated by the factor model are similar to the ones we find with UE and BB, but slightly smoother.}} 


\section{Conclusion}  \label{sec:conclusions}

UE and BB can be parametrized so that they yield the same forecasts and marginal likelihoods. However, Section~\ref{sec:FX} shows that the smoothed posteriors can be rather different even when the marginal likelihoods are identical. {{When the marginal likelihoods of the models coincide, posterior likelihood ratios, posterior predictive checks, and mixture models can be used to compare them.} }

\section*{Acknowledgement}

We thank Mike West (at Duke University) for his encouragement and feedback.

\clearpage

\appendix

\section{Distributions} \label{sec:dist}

The definitions of Wishart and matrix Beta distributions can be found, for instance, in \cite{windle} and \cite{pradowest}. We include them here for completeness. 

\textbf{Wishart:} Let $A$ be a $q \times q$ symmetric positive definite matrix. Then, $A \sim \text{Wishart}_q(h, S)$ if its probability density function  is
$$
p(A) = 2^{-{(h q)}/2} \, | S |^{-h/2} \, \Gamma_q(h/2) \, | A |^{(h - q -1)/2} \, \exp \left\{ -\frac{1}{2} \text{tr}(S^{-1} A) \right\}, 
$$
where $h > q - 1$ and $\Gamma_q(h/2)$ is the multivariate gamma function evaluated at $h/2$. The definition can be extended to $h \le q-1$, in which case $A$ is rank-deficient; see e.g. \cite{windle} for details.

\textbf{Matrix beta distribution:} Let $A_1 \sim \text{Wishart}_q(n_1, \Sigma^{-1})$ and $A_2 \sim \text{Wishart}_q(n_2, \Sigma^{-1})$ be independent random variables, where $\Sigma$ is symmetric positive-definite, $n_2 > q - 1$ and either $n_1 < q$ is an integer or $n_1 > q - 1$ is real-valued. Let $T = \text{uchol}(A_1 + A_2)$ and $B = (T^{-1})' A_1 T^{-1}$. Then, $B \sim \text{MatrixBeta}_q(n_1/2, n_2/2)$.

\section{Confluent hypergeometric function}

The confluent hypergeometric function of the second kind, $U(a,b,z)$, appears in the conditional expectation of $\Phi_{t+1}$ given $\Phi_t$. It was originally defined as the solution of specific differential equations, and there is no constraint on parameters $a$ and $b$. It admits the integral representation
\begin{equation*}
U(a,b,z) = \frac{1}{\Gamma (a)} \int _0^{\infty} e^{-zt} t^{a-1} (1+t)^{b-a-1} dt, \ \ \ \ \ \ \mathrm{for} \ a>0, \ \ \ \ \mathrm{(DLMF:13.4.4)}
\end{equation*}
In our derivations, we use Kummer's relation, 
\begin{equation*}
U(a,b,z) = z^{1-b} U(a-b+1,2-b,z) \ \ \ \ \ \ \ \mathrm{(DLMF:13.2.40)}
\end{equation*}
which is valid for any $(a,b,z)$. 

{{
\section{Derivations of main results} \label{derivations}

\subsection{Conditional expectation}

Let the Uhlig extended and beta-Bartlett models be as defined in the main text, and assume
\begin{equation*}
    k_0 = n+k,\qquad \beta = n/(n+k),\qquad b = \lambda\qquad \text{and}\qquad D^B_0 = D^U_0, \label{hp}
\end{equation*}
with $\Phi^U_{t-1} = \Phi^B_{t-1} = \Phi_{t-1} = (U_{t-1} P_{t-1})'U_{t-1} P_{t-1}$ and $U_{t-1} = (u_{ij,t-1})_{i,j\in 1:q}$. In this section, we compare $E(\Phi_t^U \mid \Phi_{t-1}, \mathcal{D}_{t-1})$ to $E(\Phi_t^B \mid \Phi_{t-1}, \mathcal{D}_{t-1})$. First, we find $E(\Phi_t^U \mid \Phi_{t-1}, \mathcal{D}_{t-1})$. Let $\Psi_t \sim \mathrm{MatrixBeta}_q(n/2,k/2)$ with expectation $E(\Psi_t) = n/(n+k)$ (see, for example, Theorem 3.2. in \cite{konno1988exact}). Then,
\begin{align*}
E(\Phi_t^U \mid \Phi_{t-1}, \mathcal{D}_{t-1}) =  (U_{t-1} P_{t-1})' \, E( \Psi_t ) \, U_{t-1} P_{t-1}/\lambda = \frac{n}{\lambda(n+k)}  P_{t-1}' U_{t-1}' U_{t-1} P_{t-1}.
\end{align*}
The derivation of $E(\Phi_t^B \mid \Phi_{t-1}, \mathcal{D}_{t-1})$ is slightly more complicated. First of all, note that
\begin{align*}
    E(\Phi_t^B \mid \Phi_{t-1}, \mathcal{D}_{t-1}) =  P_{t-1}' E( \tilde{U}_t'\tilde{U}_t \mid \Phi_{t-1}, \mathcal{D}_{t-1} ) P_{t-1}/\lambda.
\end{align*}
Now, we find the expectation of the matrix entries $E[ (\tilde{U}_t'\tilde{U}_t)_{ij} \mid \Phi_{t-1}, \mathcal{D}_{t-1} ]$, which is the only part that is missing to find $E(\Phi_t^B \mid \Phi_{t-1}, \mathcal{D}_{t-1}) $.
\begin{align*}
    E[ (\tilde{U}_t'\tilde{U}_t)_{ij} \mid \Phi_{t-1}, \mathcal{D}_{t-1} ] &= E( \sum_{l = 1}^{i \wedge j} \tilde{u}_{li,t} \tilde{u}_{lj,t} \mid \Phi_{t-1}, \mathcal{D}_{t-1} ) \\
    &= \sum_{l = 1}^{i \wedge j-1} {u}_{li,t-1} u_{lj,t-1}  + E(\tilde{u}_{i \wedge j, i \wedge j, t} \tilde{u}_{i \vee j, i \wedge j, t} ).
\end{align*}
We only need to find $ E(\tilde{u}_{i \wedge j, i \wedge j, t} \tilde{u}_{i \vee j, i \wedge j, t} )$. If $i = j$, 
$$
 E(\tilde{u}_{i \wedge j, i \wedge j, t} \tilde{u}_{i \vee j, i \wedge j, t} ) = E(\tilde{u}^2_{ii , t}  ) =  u_{ii,t-1}^2 \,  E(\eta_{i,t}) = \frac{n-i+1}{n-i+1+k} u^2_{ii,t-1}.
$$
If $i < j$,
\begin{align*}
E(\tilde{u}_{i \wedge j, i \wedge j, t} \tilde{u}_{i \wedge j, i \vee j, t} ) &= E(\tilde{u}_{ii, t} \tilde{u}_{ij, t} ) \\ 
& = u_{ij, t-1} \, E(\tilde{u}_{ii,t}) \\
& = u_{ij,t-1} u_{ii, t-1} E( \sqrt{\eta_{i,t}} ) \\
& = u_{ij,t-1} u_{ii, t-1} \frac{\Gamma((n-i+2)/2) \Gamma(n-i+k+1)/2)}{\Gamma((n-i+1)/2) \Gamma(n-i+k+2)/2)} .
\end{align*}
The expectation of $\sqrt{\eta_{i,t}}$ is straightforward to compute given that $\eta_{i,t} \sim \mathrm{Beta}((n-i+1)/2, k/2)$ (it amounts to identifying another Beta kernel). The expectation for $i > j$ can be found in an analogous manner. The expressions given in the main document are compact ways of writing the results we have derived here using minima, maxima, and Kronecker delta functions as needed.

\subsection{Example 1}

Assume Equation (3) in the main text holds. Let $k = 1$, let $P_t = \mathrm{uchol}((kD_t)^{-1})$ be the identity matrix, and $\Upsilon = \mathrm{uchol}(\mathrm{\Phi}_{t+1}) = (\upsilon_{ij})_{i,j \in 1:q}$. Let $(\Phi_t^U)_{ij} $ and $(\Phi_t^B)_{ij} $ be the the entries of $\Phi_t^U$ and $\Phi^B_t$, respectively. Then, with the Uhlig extended process:
\begin{align*}
    E[ (\Phi_t^U)_{ij}  \mid \Phi_{t+1}, \mathcal{D}_t]  &= E[ (\lambda \Upsilon'\Upsilon + Z_t)_{ij}  \mid \Phi_{t+1}, \mathcal{D}_t ] = \lambda \sum_{l=1}^{i \wedge j} \upsilon_{li} \upsilon_{lj} + \delta_{ij} \\
     V (\Phi_t^U)_{ij}  \mid \Phi_{t+1}, \mathcal{D}_t]  &= V[ (\lambda \Upsilon'\Upsilon + Z_t)_{ij}  \mid \Phi_{t+1}, \mathcal{D}_t ] = V( (Z_t)_{ij} \mid \Phi_{t+1}, \mathcal{D}_t ) = 1 + \delta_{ij},
\end{align*}
as given in the main text. 

Our next step is finding the expectation and variance for the beta-Bartlett process. For simplicity, assume $i \le j$. The case $i > j$ can be handled in an analogous manner.

Given $P_t = I$ (the identity matrix), $\Upsilon = \mathrm{uchol}(\Phi_{t+1})$, and the fact that Equation (3) holds, we have that $\tilde{U}_{t+1} = \mathrm{uchol}(\lambda \Phi_{t+1}) = \sqrt{\lambda} \Upsilon$. Consider the Bartlett decomposition $\Phi_{t} = (U_t P_t)'(U_t P_t) = U_t'U_t$. Given the state evolution described in Table 1 and the fact that $\tilde{U}_{t+1} = \sqrt{\lambda} \Upsilon$, we know that $(U_t)_{ij} = u_{ij,t}$ for $i \neq j$ is equal to $\sqrt{\lambda} \upsilon_{ij}$. For $i = j$, we have $u_{ii,t} = \sqrt{\lambda \upsilon_{ii}^2 + \theta_i}$, where $\theta_i \sim \chi^2_1$ [This result follows using standard facts for the univariate gamma-beta discount model (see e.g. Exercise 4 in Section 4.6. of \cite{pradowest}]. 

Let's examine the entries of $\Phi_t$. For $i = j$,
\begin{equation*}
\begin{split}
(\Phi _t)_{ii} &= \sum _{l=1}^i u_{li,t} u_{lj,t} \\
&= \lambda \sum _{l=1}^{i} \upsilon^2_{li} +  \theta_i,
\end{split}
\end{equation*}
so
$$
E( (\Phi _t)_{ii} \mid \Phi_{t+1}, \mathcal{D}_t) = \lambda \sum_{l=1}^i \upsilon^2_{li} + 1, \, \, \, V( (\Phi _t)_{ii} \mid \Phi_{t+1}, \mathcal{D}_t) = 2.
$$
For $i < j$:
\begin{equation*}
\begin{split}
(\Phi _t)_{ij} &= \sum _{l=1}^i u_{li,t} u_{lj,t} \\
&= \lambda \sum _{l=1}^{i-1} \upsilon_{li} \upsilon_{lj} + \sqrt{\lambda} \upsilon_{ij} \sqrt{\lambda \upsilon_{ii}^2 + \theta_i}.
\end{split}
\end{equation*}
Thus,
\begin{align*}
E[ (\Phi _t)_{ij} \mid \Phi_{t+1}, \mathcal{D}_t ] &=  \lambda \sum _{l=1}^{i-1} \upsilon_{li} \upsilon_{lj} + \sqrt{\lambda} \upsilon_{ij} E[ \sqrt{\lambda \upsilon_{ii}^2 + \theta_i} \mid \lambda, \upsilon_{ii}] \\
V[ (\Phi _t)_{ij} \mid \Phi_{t+1}, \mathcal{D}_t ] &= \lambda \upsilon_{ij}^2 \, V[ \sqrt{\lambda \upsilon_{ii}^2 + \theta_i} \mid \lambda, \upsilon_{ii}].
\end{align*}
The expectation is  
\begin{equation*}
E[ \sqrt{ \lambda \upsilon_{ii}^2 + \theta _{it} } \mid \lambda, \upsilon_{ii} ] = \int _0^{\infty} \sqrt{ \lambda \upsilon_{ii}^2 + \theta _{it} } \frac{1}{\Gamma (1/2) (1/2)^{-1/2}} \theta _{it}^{-1/2} e^{-\theta _{it}/2} d\theta _{it}
\end{equation*}
 Changing variables to $t = \theta _{it} / (\lambda \upsilon_{ii}^2)$, we obtain 
\begin{equation*}
\begin{split}
E[ \sqrt{ \lambda \upsilon_{ii}^2 + \theta _{it} }  \mid \lambda, \upsilon_{ii}] &= \int _0^{\infty} \sqrt{ \lambda \upsilon_{ii}^2}  (1+t)^{1/2} \frac{1}{\Gamma (1/2) (1/2)^{-1/2}} (\lambda \upsilon_{ii}^2)^{-1/2} t^{-1/2} e^{-\lambda \upsilon_{ii}^2t/2} (\lambda \upsilon_{ii}^2) dt \\
&= \frac{\lambda \upsilon_{ii}^2}{\sqrt{2}} \int _0^{\infty} \frac{1}{\Gamma (1/2)}  t^{-1/2} (1+t)^{1/2} e^{-\lambda \upsilon_{ii}^2t/2} dt,
\end{split}
\end{equation*}
where we read off the integral representation of $U$ with $a=1/2$, $b-a-1 = 1/2$ ($b=2$) and $z=\lambda \upsilon_{ii}^2/2$. Thus, using Kummer's formula, we have 
\begin{equation*}
\begin{split}
E[ \sqrt{ \lambda \upsilon_{ii}^2 + \theta _{it} }  \mid \lambda, \upsilon_{ii}] &= \frac{\lambda \upsilon_{ii}^2}{\sqrt{2}} U(1/2,2,\lambda \upsilon_{ii}^2/2) \\
&= \frac{\lambda \upsilon_{ii}^2}{\sqrt{2}} (\lambda \upsilon_{ii}^2/2)^{1-2} U(1/2-2+1,2-2,\lambda \upsilon_{ii}^2/2) \\
&= \sqrt{2} U(-1/2,0,\lambda \upsilon_{ii}^2/2) 
\end{split}
\end{equation*}
In summary, we obtain 
\begin{equation*}
\begin{split}
E[ (\Phi _t)_{ij} \mid \Phi_{t+1}, \mathcal{D}_t ] 
&= \lambda \sum _{l=1}^{i-1} \upsilon_{li}\upsilon_{lj} + \sqrt{ 2\lambda }\upsilon_{ij} U(-1/2,0,\lambda \upsilon_{ii}^2/2), 
\end{split}
\end{equation*}
as shown in the main text. The variance can be computed similarly. We have 
\begin{equation*}
\begin{split}
V[ (\Phi _t)_{ij} \mid \Phi_{t+1}, \mathcal{D}_{t+1} ] &= \lambda \upsilon_{ij}^2  V[ \sqrt{ \lambda \upsilon_{ii}^2 + \theta _{it} } \mid \lambda, \upsilon_{ii} ] \\
&= \lambda \upsilon_{ij}^2 \{ \ E[ \lambda \upsilon_{ii}^2 + \theta _{it} \mid \lambda, \upsilon_{ii} ] - E[ \sqrt{ \lambda \upsilon_{ii}^2 + \theta _{it} } \mid \lambda, \upsilon_{ii} ]^2 \  \} \\
&= \lambda \upsilon_{ij}^2 \{ \ \lambda \upsilon_{ii}^2 + 1 - E[ \sqrt{ \lambda \upsilon_{ii}^2 + \theta _{it} } \mid \lambda, \upsilon_{ii} ]^2 \  \}.
\end{split}
\end{equation*}
The expressions given in the main text can be found by putting together the terms derived here, using maxima and minima as needed.
}}
\section{Additional example} \label{addexample}

 Let $\Phi_{t+1} = \text{diag}(\phi_{1:q})$ and $D_{t}^{-1} = \text{diag}(d_{1:q})$. For simplicity, we set $k = 1$, although the same computations could be done for $k \neq 1$. For the UE process,  
\begin{align*}
 E( \Phi^U_t \mid \Phi_{t+1}, D_t) = \text{diag}(\lambda \phi_{1:q}+d_{1:q}); \, \, \, V[ (\Phi^U_{t})_{ij} \mid \Phi_{t+1}, D_t]  = \delta_{ij} \, d_i^2 + d_i d_j,
\end{align*}
where $\delta_{ij}$ is Kronecker's delta function. For the BB process, we have
 \begin{align*}
 E( \Phi^B_t \mid \Phi_{t+1}, D_t) &= \text{diag}(\lambda \phi_{1:q}+d_{1:q}); \, \, \, V[ (\Phi_t)^B_{ij} \mid \Phi_{t+1}, D_t] = \delta_{ij} \, 2 d_i^2.
\end{align*}
 The conditionals are equal in expectation, but the variance of the off-diagonal elements don't coincide.  With the BB process, the off-diagonal elements are 0 with probability 1, whereas with the UE process the off-diagonal elements aren't identically equal to 0.

\section{Foreign exchange rates application} \label{appFX}

\subsection{Technical details}

As we mentioned in the main text, we take $k = 1$, which implies that we can simply work with normal likelihoods for the returns.

In \cite{windle}, the discounting parameter $\lambda$ is automatically chosen to satisfy $\lambda ^{-1} = 1 + k / (n-q-1)$. This constraint not only reduces the number of parameters to estimate, but also guarantees that $E[\Phi _t^{-1}|\mathcal{D}_t] = E[\Phi _{t+1}^{-1}|\mathcal{D}_t]$, a property the authors deem desirable. In contrast, we directly maximize the marginal likelihood with respect to $(n,\lambda)$ on a grid under no constraint. The marginal likelihood is the product of one-step ahead forecast densities, each of which is the multivariate-$t$ distribution defined by
\begin{equation*}
\begin{split}
p (r_t |\mathcal{D} _{t-1} ) &= \int N_q(r_t|0,\Phi _t^{-1}) W_q(\Phi _t | n, (\lambda D_{t-1})^{-1}) (d\Phi _t) \\
&= \frac{\Gamma ( n/2 )}{ \Gamma ( (n+1-q)/2 ) } \frac{ |\lambda D_{t-1}|^{-1/2}}{\pi ^{q/2}} (1+r_t'D_{t-1}^{-1}r_t / \lambda )^{ -(n+1)/2 }. 
\end{split}
\end{equation*}
In addition, the determinant of $D_t$ is sequentially updated using the convenient relation
\begin{equation*}
    \log|D_t| =  \log (1+r_t'D_{t-1}^{-1}r_t / \lambda ) + q\log (\lambda ) + \log | D_{t-1}| ,
\end{equation*}
so the evaluation of marginal likelihood isn't computationally demanding. For maximizing the marginal likelihood, we evaluate it at $n\in \{ 3,4,\dots , 20 \}$ and $\lambda \in \{ 0.600, 0.601, \dots , 0.990 \}$. 

The posterior likelihood ratio \citep{aitkin1991posterior} can be hard to estimate numerically, but its logarithm is stable. To see this, recall that
$$
\overline{\ell}_{UB} = \log\left\{E[ L(\Phi^U_{0:T}) | \mathcal{D}_T]/E[ L(\Phi^B_{0:T}) | \mathcal{D}_T] \right\}, \, \, \, L(\Phi_{0:T}) = \prod_{t=1}^T N_q(r_t \mid 0_q, \Phi_t^{-1}).
$$
Based on Monte Carlo samples $\Phi^{U\ast}_{1:N} \sim \Phi^U_{0:T} \mid \mathcal{D}_T$ and $\Phi^{B\ast}_{1:N} \sim \Phi^B_{0:T} \mid \mathcal{D}_T$,
$$
\overline{\ell}_{UB} \approx \text{LSE}( \ell(\Phi^{U\ast}_{1:N})) - \text{LSE}( \ell(\Phi^{B\ast}_{1:N})),
$$
where $\ell$ is $\log L(\Phi_{0:T})$ and LSE is the log-sum-exp function, which can be implemented in a numerically stable way. 

We implement the mixture model approach proposed in \cite{kamary2014testing} through a missing-data augmented Gibbs sampler (see Appendix~\ref{appFX}). The target model is defined by the mixture of likelihoods, 
\begin{equation*}
r_t  \mid \alpha , \Phi _{t}^U, \Phi _{t}^B \ \sim \ \alpha  N_q(0,(\Phi _t^U)^{-1} ) + (1-\alpha ) N_q(0,(\Phi_t^B)^{-1} ).
\end{equation*}
We implement the following augmented model: 
\begin{equation*}
\begin{split}
&r_t = z_i r_t^U + (1-z_i) r_t^B \\
&r_t^{\mathcal{M}} \mid \Phi _t^{\mathcal{M}} \sim N_q(0,(\Phi _t^{\mathcal{M}})^{-1}), \ \ \ \ \ \ \mathcal{M} \in \{ U, B \} \\
&z_i \stackrel{\text{iid}}{\sim} \text{Bernoulli}(\alpha ) \\
&\alpha \sim \text{Beta}(a_0,b_0)
\end{split}
\end{equation*}
The actual observed return, $r_t$, is defined separately from the inputs of two models, $r_t^U$ and $r_t^B$. At each iteration of the Gibbs sampler, conditional on $z_i$, we decide which model is fed by $r_t$, and which model is ``missing'' its observation. The notable advantage of this approach is that the missing observation, either $r_t^U$ or $r_t^B$, is a parameter, so it is sampled through the course of the Gibbs sampler. As a result, the sampling of $\Phi _{1:T}^{\mathcal{M}}$ is based on a full sequence of observations $r_{1:T}^{\mathcal{M}}$  and we can apply the forward filtering equations and backward sampler we described in the main text.

The Gibbs sampler consists in iteratively sampling from the following full-conditional distributions:
\begin{itemize}
	\item Sample $z_t$ from Bernoulli distribution with 
	probability 
	\begin{equation*}
	\begin{split}
	P[ z_t=1 | - ] &\propto \alpha N_q(r_t^U|0,(\Phi _t^U)^{-1} ) \\
	P[ z_t=0 | - ] &\propto (1-\alpha ) N_q(r_t^B|0,(\Phi _t^B)^{-1} ).
	\end{split}
	\end{equation*}
	 It's more computationally stable to work on the log-scale:
	\begin{equation*}
	\begin{split}
	\log P( z_t=1 | - ) &= c + \log( \alpha ) + \frac{1}{2} \log | \Phi _t^U | - \frac{1}{2} (r_t^U)'\Phi _t^U r_t^U \\
	\log P( z_t=0 | - ) &= c + \log( 1 - \alpha ) + \frac{1}{2} \log | \Phi _t^B | - \frac{1}{2} (r_t^B)'\Phi _t^B r_t^B,
	\end{split}
	\end{equation*}
	where $c$ is the common constant. 
	
	\item Define $r_t^U$ and $r_t^B$ as follows: 
	\begin{itemize}
		\item If $z_t = 1$, then set $r_t^U = r_t$ and generate $r_t^B \sim N_q(0,(\Phi _t^B)^{-1})$. 
		\item If $z_t = 0$, then generate $r_t^U \sim N_q(0,(\Phi _t^U)^{-1})$ and set $r_t^U = r_t$.
	\end{itemize} 
	
	\item Sample $\alpha$ from $\text{Beta}( a_1 , b_1 )$, 
	\begin{equation*}
	a_1 = a_0 + \sum _{t=1}^T z_t, \ \ \ \ \ \ 	b_1 = b_0 + \sum _{t=1}^T (1-z_t)
	\end{equation*} 
	
	\item Sample $\{ \Phi _{1:T}^U \}$ and $\{ \Phi _{1:T}^B \}$ using the forward-filtering equations and the backward sampler described in the main text. 
	
\end{itemize}

\subsection{Additional results and figures}

Figure~\ref{fig:returns} shows the original series of returns. Figure~\ref{fig:ml} shows the contours of the marginal likelihood, along with the maximizer $(n,\lambda )=(5,0.799)$ indicated by the red circle. In this figure, we also show the maximizer under the constraint that was used in \cite{windle}, $(n,\lambda ) = (10,0.857)$, indicated by the blue box. The posterior and predictive analysis in this study is based on the former choice.

\begin{figure}[h!]
	\centering
	\includegraphics[width=4.5in]{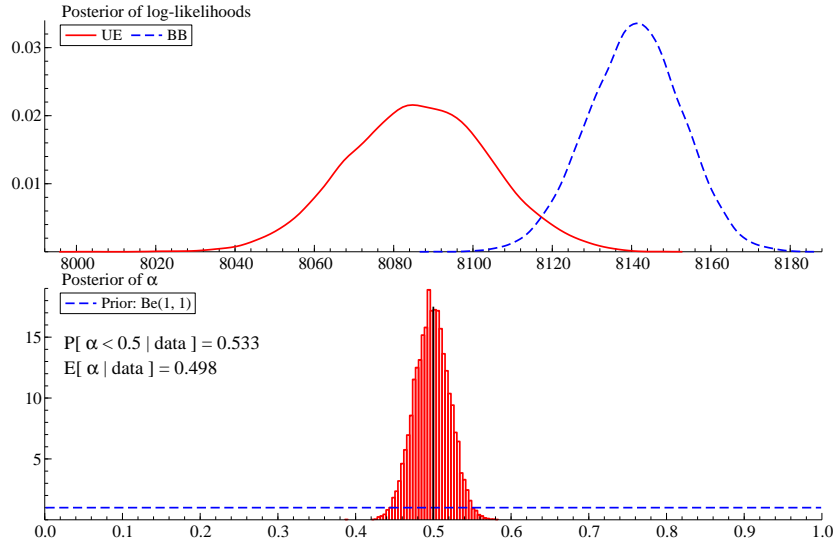}
	\caption{ Top: Posterior log-likelihoods with Uhlig extended (UE; solid red) and beta-Bartlett (BB; dashed blue). Bottom: Prior on mixture weight $\alpha$ (dashed blue) and posterior (red histogram). Vertical line at 0.5.
	} \label{fig:plrmixt}
\end{figure}

Figure~\ref{fig:plrmixt} visualizes the results with posterior likelihood ratios and the mixture model approach of \cite{kamary2014testing}. The top plot shows that the posterior log-likelihood of the beta-Bartlett model is clearly higher than that of the Uhlig extended model. The bottom plot shows that the mixture approach does prefer the beta-Bartlett model as well, but to a much lesser extent. 

\begin{figure}[h!]
	\centering
	\includegraphics[width=4.5in]{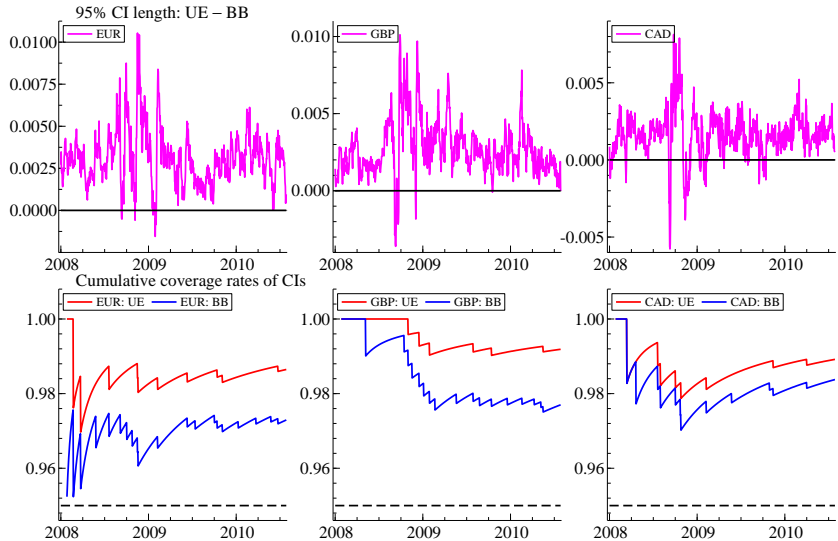}
	\caption{ Length of 95\% posterior predictive intervals and cumulative empirical coverage rates. 
	} \label{fig:range}
\end{figure}

\begin{figure}[h!]
	\centering
	\includegraphics[width=4.5in]{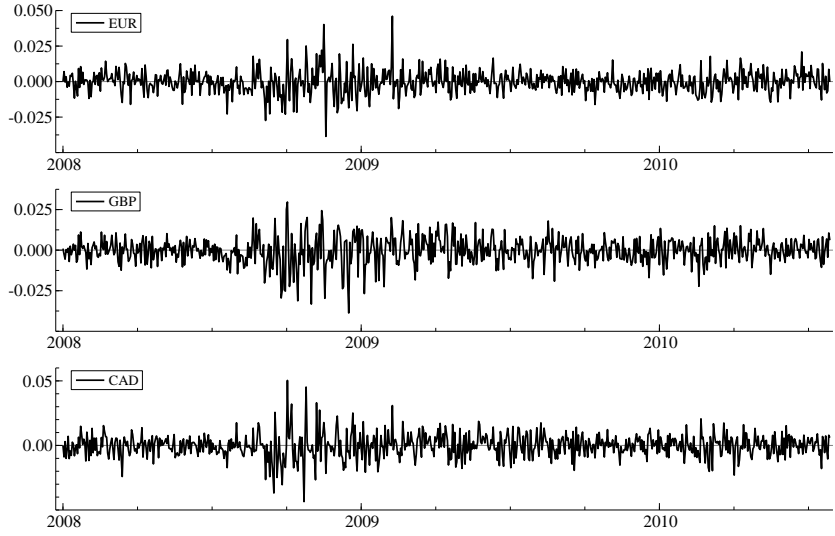}
	\caption{Time series of daily returns from EUR, GBP and CAD in US dollars.
	} \label{fig:returns}
\end{figure}

\begin{figure}[h!]
	\centering
	\includegraphics[width=4.5in]{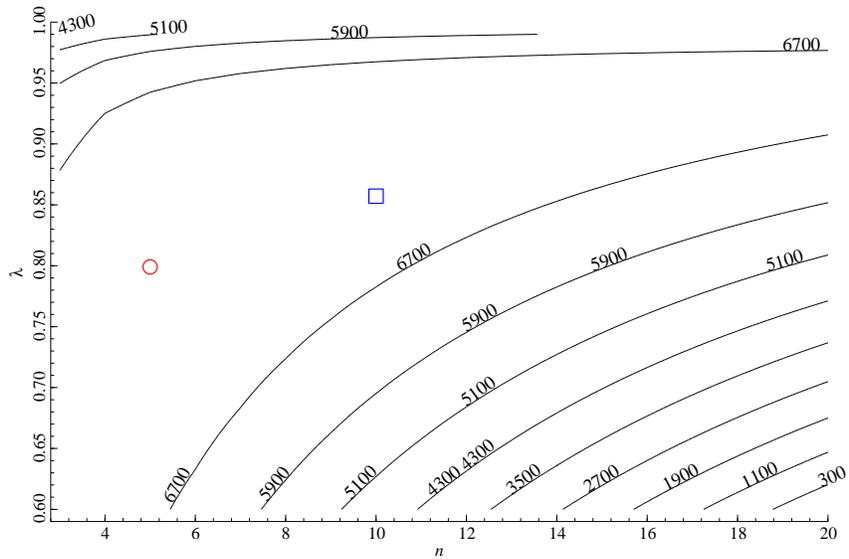}
	\caption{ Contour plots of marginal likelihoods as the functions of $(n,\lambda )$. The red circle indicates the maximizer $(5,0.799)$.  The blue box shows $(10,0.857)$, the maximizer under the constraint. 
	}  \label{fig:ml}
\end{figure}

Figure~\ref{fig:range} shows lengths of 95\% posterior predictive intervals and their cumulative empirical coverage over time. The beta-Bartlett process reports shorter intervals most of the time, and their empirical coverage is always above 95\%.

\begin{figure}[h!]
	\centering
	\includegraphics[width=4.5in]{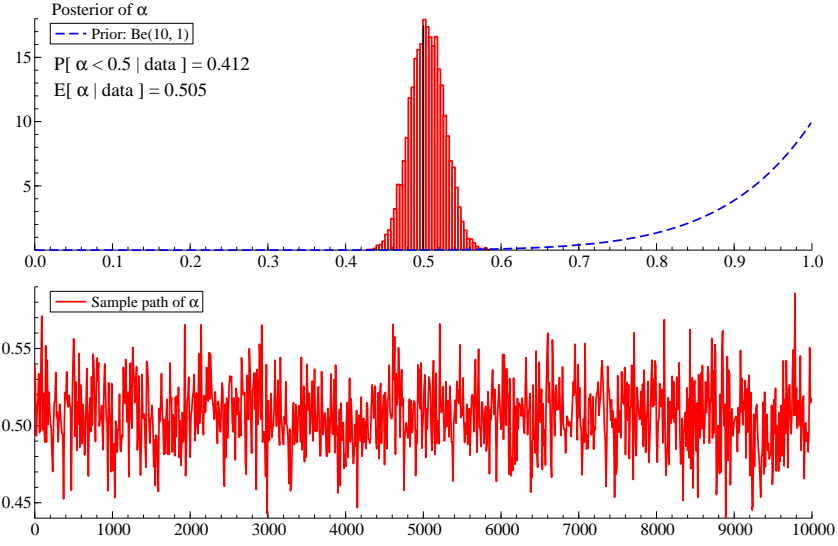}
	\caption{Histogram and sample path of $\alpha$ for prior $\alpha \sim Be(10,1)$
	} \label{fig:weight2} 
\end{figure}

 In the main text, we presented our results with the mixture model approach in \cite{kamary2014testing} with $\alpha \sim \text{Beta}(1,1)$. We also tried other hyperparameters to test out the effect of the prior.  For example, consider $\alpha \sim \text{Beta}(10,1)$, so the UE model is strongly preferred a priori. The results  are shown in Figure~\ref{fig:weight2}. The estimated posterior mean is $0.505123$ (standard error: $0.000261316$). The estimated posterior probability of having $\alpha < 0.5$ is $0.4124$ (standard error: $0.00537318$). This shows that the prior concentrates towards roughly 0.5 even if the starting point is far away from it. 

\begin{figure}[!htbp]
	\centering
	\includegraphics[width=5in]{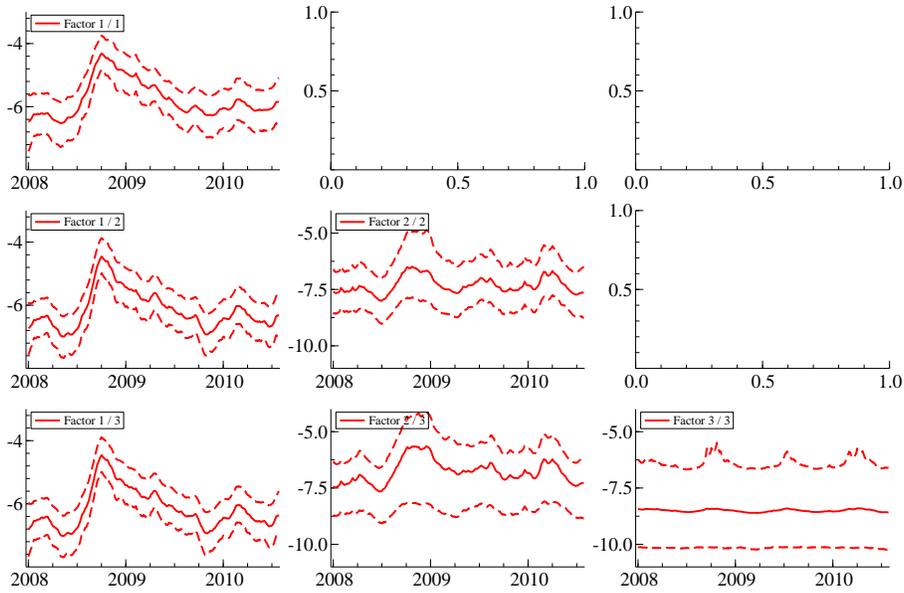}
	\caption{Factor log-volatilities $\lambda _{it}$ under  $k=1$ (top row), $k=2$ (middle row) and $k=3$ (bottom row).}
	\label{fig:vol}
\end{figure}

\begin{figure}[!htbp]
	\centering
	\includegraphics[width=5in]{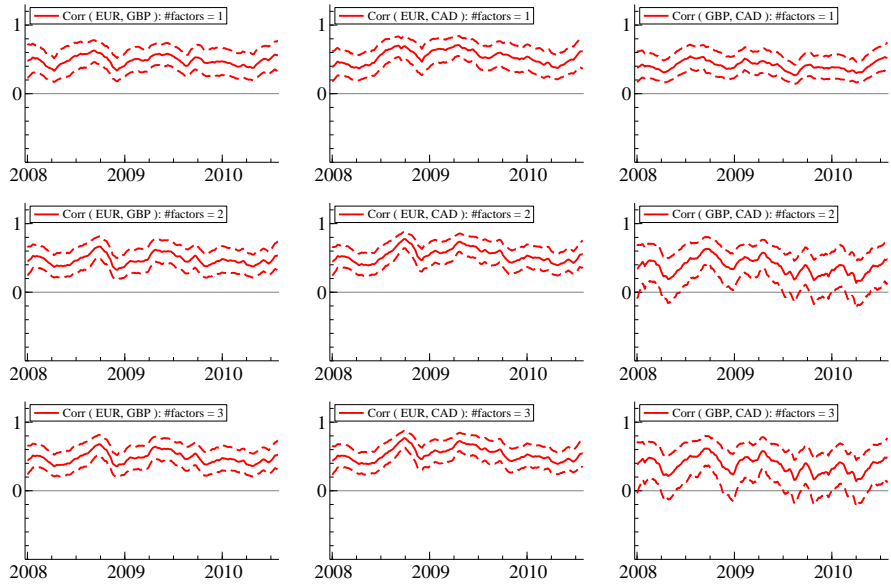}
	\caption{Posteriors of the pairwise correlations of returns, computed by the sampled covariance matrix $\Sigma _t = XH_tX' + \Psi _t$ under models with $k=1$ (top row), $k=2$ (middle row) and $k=3$ (bottom row).}
	\label{fig:cor}
\end{figure}

{{\subsection{Comparison to a Bayesian factor model}

In this section, we fit a multivariate factor stochastic volatility on the FX data and compare it to the results we find with UE/BB models. Our factor model is based on \cite{aguilar2000bayesian}. Extensions of this model can be considered (for example, we can add dynamic loadings as in \cite{lopes2007factor}), but we believe that this factor model is sufficiently flexible for our study. 

For $q$-dimensional returns $y_t = (y_{1t}, \, \dots \, , y_{qt})'$, the sampling model is specified as 
\begin{equation*}
y_t = X_t f_t + \epsilon _t, \ \ \ \ \ \mathrm{where} \ \ \ \ \ \epsilon _t\sim N(0_q,\Psi _t), 
\end{equation*}
where $f_t$ is the vector of $k$ factors ($k{\times}$1), $X_t$ is the factor loading matrix ($q{\times}k$) and $\Psi _t = \mathrm{diag} ( \psi _{1t},\dots , \psi _{qt} )$. We consider a constant loading matrix, i.e., $X_t = X$ for all $t$. For identifiability, the loading matrix is assumed to be lower-triangular with diagonal unity. In our application, where $y_t$ is the returns from the three FX rates and $q=3$, the loading matrix $X$ is specified as 
\begin{equation*}
X = \begin{bmatrix}
1 \\
x_{21} \\
x_{31} 
\end{bmatrix} \ (k=1), \ \ \ \ \ \begin{bmatrix}
1 & 0 \\
x_{21} & 1 \\
x_{31} & x_{32}
\end{bmatrix} \ (k=2), \ \ \ \ \ \begin{bmatrix}
1 & 0 & 0 \\
x_{21} & 1 & 0 \\
x_{31} & x_{32} & 1 
\end{bmatrix} \ (k=3),
\end{equation*}
depending on the number of factors $k$ used in the model. 
The factors are conditionally mutually independent and Gaussian, 
\begin{equation*}
f_t \stackrel{\mathrm{iid}}{\sim} N(0_k,H_t), \ \ \ \ \ \mathrm{where} \ \ \ \ \ H_t = \mathrm{diag} (h_{1t},\dots , h_{kt}).
\end{equation*}
The variance parameters in the diagonal entries of $H_t$ and $\Psi _t$ are modeled independently with univariate stochastic volatility models \citep{jacquier1994vbayesian}. That is, for $\lambda _{it} = \log \ h_{it}$, we assume
\begin{equation*}
\lambda _{it} = (1-\phi _i ) \mu _i + \phi _i \lambda _{it} + N(0,\sigma ^2_i),
\end{equation*}
for each $i=1, \, \dots \, , k$, where the triplet of AR(1) parameters, $(\mu _i,\phi _i, \sigma ^2_i)$, follows a prior which we specify later in the document. The initial value of the log-volatilities follows its stationary marginal, i.e., $h_{i0} \sim N(\mu _i,\sigma _i^2 /(1-\phi^2_i))$. Likewise, for $\eta _{it} = \log \ \psi _{it}$, we have 
\begin{equation*}
\eta _{it} = (1-\rho _i ) \alpha _i + \rho _i \eta _{it} + N(0,\xi ^2_i),
\end{equation*}
for $i=1, \, \dots \, , q$, where the AR(1) parameters are denoted by $(\alpha _i, \rho _i, \xi ^2_i)$. 

A Markov Chain Monte Carlo (MCMC) sampling method can be easily built for this model. For sampling the log-volatilities, we use a multi-move sampler \citep{shephard1997likelihood,watanabe2004multi}, as implemented in \cite{nakajima2013bayesian}. The other parameters are conditionally conjugate and easily sampled from their full conditionals directly (or by utilizing Metropolis-Hasting steps, especially when sampling $\phi_i$). The computation is implemented in Ox, based on the code used in \cite{nakajima2013bayesian} that is publicly available. 

We use proper priors for all the parameters in the model above. The hyperparameters are mostly based on choices made in  \cite{nakajima2013bayesian}. 
\begin{itemize}
	\item AR(1) coefficient $\phi _i$ (and $\rho _i$).
	\begin{equation*}
	\frac{\phi _i+1}{2} \stackrel{\mathrm{iid}}{\sim} \mathrm{Beta}(20,1.5).
	\end{equation*}
	
	\item AR(1) location $\mu _i$ (and $\alpha _i$). In the exponential scale, where $\gamma _i = \exp \{ \mu _i \}$, we set 
	\begin{equation*}
	\gamma_i \stackrel{\mathrm{iid}}{\sim} \mathrm{Gamma}(3,0.03),
	\end{equation*}
	where $\mathrm{Gamma}(a,b)$ is the gamma distribution with shape $a$ and rate $b$ (mean $a/b$).
	
	\item AR(1) variance $\sigma ^2_i$ (and $\xi ^2 _i$)
	\begin{equation*}
	\sigma ^{-2}_i \stackrel{\mathrm{iid}}{\sim} \mathrm{Gamma}(2,0.01)/
	\end{equation*}
	
	\item Factor loadings $x_{ij}$ for $i>j\ge 2$: 
	\begin{equation*}
	x_{ij} \stackrel{\mathrm{iid}}{\sim} N(0,1).
	\end{equation*}
	
\end{itemize}

We obtained 20,000 posterior samples after 2,000 of burn-in. It took about 8 minutes on a workstation to complete the posterior sampling under the factor model with $k=3$, while it took less than a second under the UE/BB models on a laptop computer. 

Figure~\ref{fig:vol} shows the trajectories of posterior means and 95\% credible intervals of factor log-volatilities $\lambda _{it}$ under the three factor models with $k=1,2,3$. In the model with $k=3$ factors, the third log-volatilities are almost constant and very negatively-large, implying that the corresponding factor is very close to zero  and, hence, ignorable. The second factor is slightly higher than the third one, but substantially smaller than the first one and less volatile.  We conclude that the second and third factors are ``redundant'' in this application, which makes sense in context, since we are working with highly-related Western currencies. 

Figure~\ref{fig:cor} summarizes the posterior correlations under the three factor models. The plots are all similar, except for the correlation between GBP and CAD where we observe more uncertainty as we have more factors in the model. The similarity of posterior analyses under the three models supports the use of the model with one factor ($k=1$) for its simplicity. In comparison with the posterior under the UE and BB models (Figure~1 in the main text), we find that the pattern of correlation dynamics resembles those of the UE and BB models, but the posteriors under the factor model are less dynamic than those of the UE and BB models. 

 Now, we compare UE/BB and the factor model in terms of predictive uncertainty quantification. The point forecasts by those models are always zero by definition, so it is not meaningful to look at the predictive MSEs. We computed the 90\% and 95\% credible intervals of one-step ahead predictive distributions, i.e., $p(y_{i,t+1}| \mathcal{D}_t)$, under the factor model ($k=1$) and the two UE/BB models $(n=5,8)$ in 64 days between January and March in 2010. The prediction by the factor model relies heavily on the repeated use of the MCMC samples. In our example, which is of very moderate dimension, it took hours to run. In contrast, the predictive marginal under the UE/BB models have a closed-form expression: $y_t|\mathcal{D}_t$ follows a multivariate $t$-distribution, and $y_{it}|\mathcal{D}_t$ are univariate $t$-distributions. 

In Figure~\ref{fig:pred}, the predictive intervals are plotted along with the actual observations. The dynamics of predictive intervals under the factor model, plotted in the left column, are almost identical due to having a single factor shared among the three currencies. In terms of empirical coverage, the predictive uncertainty expressed by these intervals is slightly overestimated, especially in the first half period in predicting GBP. The prediction would be more flexible and more dynamic if we increased the number of factors, at the cost of inflated predictive uncertainty. We also see that the dynamic intervals are less smoothed over time, due to Monte Carlo errors. 

UE/BB models with $n=5$, as used in the main text, clearly overestimates the predictive uncertainty. This result is not surprising; we chose the hyperparameters for this model by maximizing the marginal likelihood using the whole period of observations (2008-2012). That is, the model is optimized to explain the entire time series, including the extremely high volatility in 2008, which explains the inflated predictive uncertainty after 2008. By showing the analysis by another choice of hyperparameter ($n=8$), we demonstrate that one can calibrate the predictive uncertainty to an appropriate level. Overall, we see an advantage (disadvantage) of UE/BB models in its flexibility (myopic adaptation) to sample variations. This suggests that UE/BB models are mostly appropriate for short-term prediction. }}

\clearpage

\begin{figure}[!htbp]
	\centering
	\includegraphics[width=5in]{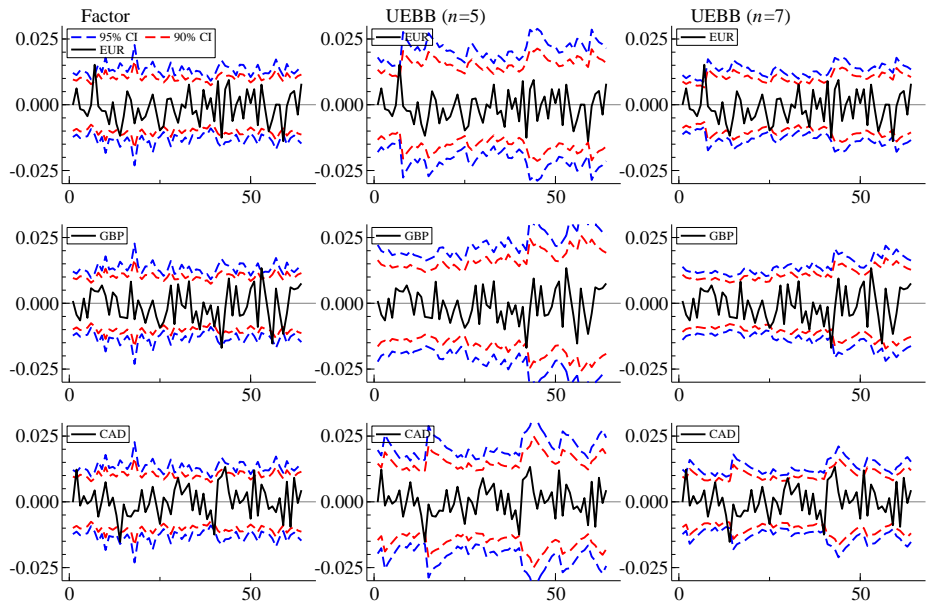}
	\caption{Predictive 90\% (red) and 95\% (blue) credible intervals of the 64 returns from the three currencies in Jan-Mar 2010 under the factor stochastic volatility models (left column), the UE/BB model with the choice of hyperparameter $n=5$ (middle column), and $n=8$ (right column).}
	\label{fig:pred}
\end{figure}%

\clearpage

\bibliographystyle{chicago}
\bibliography{main}

\end{document}